\documentclass[11pt]{amsart}
\usepackage{amssymb,amsmath}
\usepackage{amsthm}
\usepackage{amsfonts}
\usepackage[dvips]{graphicx}
\usepackage{enumerate}

\newcommand{\T}{\mathbb{T}}
\newcommand{\R}{\mathbb{R}}
\newcommand{\Z}{\mathbb{Z}}
\newcommand{\C}{\mathbb{C}}

\renewcommand{\epsilon}{\varepsilon}

\newcommand{\GNZ}{\mbox{GL}(N,\Z)}
\newcommand{\GR}[1]{\mbox{GL}(#1,\R)}
\newcommand{\GZ}[1]{\mbox{GL}(#1,\Z)}
\newcommand{\Diff}[1]{\mathit{Diff}(#1)}

\newtheorem{theorem}{Theorem}[section]

\newtheorem{proposition}[theorem]{Proposition}
\newtheorem{corollary}[theorem]{Corollary}

\newtheorem{lema}[theorem]{Lemma}

\newtheorem{remark}{Remark}

\newtheorem{definition}[theorem]{Definition}

\newtheorem{problem}{Problem}

\title[Global rigidity of certain abelian actions.]{Global rigidity of certain abelian actions by toral automorphisms.}
\author{Federico Rodriguez Hertz}
\thanks{This work was partially supported by FCE 9021, CONICYT-PDT 29/220 and CONICYT-PDT 54/18 grants}
\address{IMERL, Montevideo, Uruguay}
\email{frhertz@fing.edu.uy}
\date{\today}
\begin{document}

\begin{abstract}
We prove global rigidity results for some linear abelian actions on
tori. The type of actions we deal with includes in particular
maximal rank semisimple actions on $\T^N$.
\end{abstract}

\maketitle
\begin{section}{Introduction}
Let $\Gamma$ be a subgroup of $\mbox{GL}(N,\Z)$, the group of
$N\times N$ matrices with integral entries and determinant $\pm 1$.
We can see $\Gamma$ as acting on $\T^N=\R^N/\Z^N$ by matrix
multiplication. In this case we will say that $\Gamma$ induce a
linear action or the standard action on $\T^N$. In general, an
action of $\Gamma$ on $\T^N$ will be an embedding
$\rho:\Gamma\to\Diff{\T^N}$ and we will say that $\rho$ is an Anosov
action if it has an Anosov element, i.e. if there is $m\in\Gamma$
such that $\rho(m)$ is an Anosov diffeomorphism. In this paper we
will be concerned with global rigidity results for abelian linear
actions on $\T^N$. We shall say that the standard action of $\Gamma$
on $\T^N$ is globally rigid if any Anosov action of $\Gamma$ on
$\T^N$ which induces the standard action in homology is smoothly
conjugated to it.

\begin{theorem}\label{cent}
Let $A\in\GNZ$, be a matrix with characteristic polynomial
irreducible over $\Z$. Assume also that the centralizer $Z(A)$ of
$A$ in $\mbox{GL}(N,\Z)$ has rank at least $2$. Then the associated
action of any finite index subgroup of $Z(A)$ on $\T^N$ is globally
rigid.
\end{theorem}

We want to remark that due to the Dirichlet unit theorem, in the
above case, $Z(A)$ is a finite extension of $\Z^{r+c-1}$ where $r$
is the number of real eigenvalues and $c$ is the number of pairs of
complex eigenvalues, $r+2c=N$. So, $Z(A)$ has rank one only if $N=2$
or if $N=3$ and $A$ has a complex eigenvalue or if $N=4$ and $A$ has
only complex eigenvalues.

We think that the following should also be true.
\begin{problem}
Let $\Gamma$ be any finite index subgroup of $Z(A)$ for
$A\in\mbox{GL}(N,\Z)$, $N\geq 3$, assume also that $Z(A)$ is big
enough. Is the standard action of $\Gamma$ on $\T^N$ globally rigid?
\end{problem}

Observe that when $A$ is the identity matrix, $\Gamma$ is any finite
index subgroup of $\mbox{GL}(N,\Z)$, see \cite{klz}. Also, one may
formulate the local rigidity problem and a similar problem for
actions on infra-nilmanifolds.

Classification of Anosov actions is one of the most striking
problems in dynamics. When $\Gamma=\Z$, thus the action is generated
by a diffeomorphism $f$, Franks \cite{fr} and Manning \cite{ma} have
proven that if $M$ is a torus, a nilmanifold or an infra-nilmanifold
then $f$ is topologically conjugated to an automorphism and thus it
is essentially of an algebraic nature. In \cite{ne}, Newhouse proved
that codimension one Anosov diffeomorphism always live in tori. On
the other hand, Brin in \cite{br} get that with some bunching
hypothesis in the spectrum of the differential of $f$ the manifold
should be an infra-nilmanifold also. It is conjectured that Anosov
diffeomorphisms are always of algebraic nature, up to topological
conjugacy.

When dealing with higher rank actions typically more can be said,
for example that the topological conjugacy is smooth. At least this
is true when $\rho$ is a small perturbation of an "irreducible"
algebraic Anosov action of $\Z^k$, $k\geq 2$, see Katok and Spatzier
\cite{ks}. Moreover, it is conjectured \cite{kas}, that every
"irreducible" $\Z^k$, $k\geq 2$, Anosov action on any compact
manifold is smoothly conjugated to an algebraic action.

Theorem~\ref{cent} has an interesting particular case that is when
dealing with Cartan actions, that is, when the matrix $A$ has only
real eigenvalues. This case was already studied by Katok and Lewis
in \cite{kl} where the local rigidity property and also some global
rigidity but with some restriction on the nonlinear action was
established. In general one can try to push this notion of Cartan
action into a broader non linear context asking the manifold $M$ to
splits into $d$ invariant directions, $d=dim (M)$ and with the
action having Anosov elements that contracts and expands this
directions. Some global rigidity properties for this type of actions
were studied recently by Kalinin and Spatzier in \cite{kas}.

On the other hand, one can study the measure rigidity problem and in
this case, already for the linear action it is not known if the
unique invariant measures are Lebesgue and the atomic ones. This was
first noticed by H. Furstenberg \cite{fu} who posed the problem of
whether the unique invariant measures for the $\times 2$ $\times 3$
action on the circle are Lebesgue and the atomic ones. What is known
is that when the entropy of the measure is positive then the measure
should be Lebesgue, see \cite{ru}, \cite{jo} for the $\times 2$
$\times 3$ case and \cite{ks2} for the case of Cartan actions. There
is also lot of work on the study of the measure rigidity for linear
actions, see \cite{kak1} and \cite{kak} for a good account an
references about this case. But for the nonlinear case there is not
much work, in \cite{kak} Kalinin and Katok remarkably have proven
that $\Z^k$ actions on $\T^{k+1}$ (a priori not Anosov actions) that
induce a Cartan action in homology should leave invariant a measure
absolutely continuous w.r.t. Lebesgue. Later, Katok with the author
in \cite{krh} proved that this measure is unique in some sense and
that the action is in fact measurably isomorphic to the linear one.
For the general nonlinear case, Kalinin, Katok and the author
\cite{kakrh} prove the existence of an invariant measure absolutely
continuous w.r.t. Lebesgue for quite general actions of $\Z^k$ on a
$k+1$ dimensional manifold.

This paper grew up from a conversation with Anatole Katok during a
visit of the author to the Penn State University in October 2001. At
that time he told me the problem of the global rigidity of $\Z^2$
actions on $\T^3$, this case was solved in \cite{rh}. The big step
from there to here is to get rid of the one dimensionality of the
invariant spaces.

I would like to thank Anatole Katok for introducing me into the
wonderful subject of rigidity and also the people at Penn State for
their warm hospitality.

In section~\ref{basics} we shall expose some definitions that will
be needed for the paper and state the main theorem. We recommend the
reader to jump to the last section to see how is the scheme of the
proof of the main theorem before reading sections~\ref{smolin} and
\ref{righol}. Finally section~\ref{varios} is not about higher rank
actions and applies to single diffeomorphisms. We think that each of
sections~\ref{varios},~\ref{smolin} and~\ref{righol} have their
independent interest and they are in fact quite independent.
\end{section}




\begin{section}{Theorems and definitions}\label{basics}
In this section we shall expose the basics notions and state the
main theorem.
\begin{subsection}{Definitions}
\begin{subsubsection}{Conjugacies}
Let $f:\T^N\to\T^N$ be an Anosov diffeomorphism and let $A\in\GNZ$
be its action in homology. By the results in \cite{fr} and \cite{ma}
there is a unique conjugacy $h:\T^N\to\T^N$, $h\circ f=Ah$,
homotopic to the identity, moreover, $h$ and $h^{-1}$ are H\"older
continuous. If $\rho:\Gamma\to\Diff{\T^N}$ is an abelian action with
an Anosov element then the above mentioned conjugacy will conjugate
the whole action, that is, if $\rho_*:\Gamma\to\GNZ$ is the induced
action in homology then $h\circ\rho=\rho_* h$. So that to prove the
rigidity results we shall see that $h$ is a diffeomorphism.

Given an abelian action $\rho:\Gamma\to\Diff{\T^N}$ and a point
$p\in\T^N$, let $\Gamma_p$ be the \emph{stabilizer of $p$}, that is
$\Gamma_p=\{n\in\Gamma\;:\;\rho(n)(p)=p\}$. When $\Gamma_p$ is a
finite index subgroup of $\Gamma$ we say that $p$ is a
\emph{periodic point} and call $\Gamma/\Gamma_p$ its \emph{period}.
When $\rho$ is an abelian Anosov action the periodic points for the
action coincide with the periodic points for the Anosov element.

As we can always take a finite index subgroup of $\Gamma$ isomorphic
to $\Z^k$ and we have to prove that the conjugacy $h$ is
differentiable, we will be working typically with $\Z^k$ actions.
Also, we will deal indistinctly with an action
$\rho:\Gamma\to\Diff{\T^N}$ or its image subgroup
$\Gamma\sim\rho(\Gamma)\subset\Diff{\T^N}$ and when working with a
linear action we will simply denote $\rho_*:\Gamma\to\GR{d}$ or its
image subgroup $\Gamma\sim\rho(\Gamma)\subset\GR{d}$, idem for
$\GZ{N}$.

\end{subsubsection}

\begin{subsubsection}{Lyapunov exponents}
Let $\Gamma\subset\GR{d}$ be a subgroup isomorphic to $\Z^k$. Let us
use the letter $\chi$ to denote the Lyapunov exponents of the action
induced by $\rho:\Z^k\to\Gamma\subset\GR{d}$, hence $\chi_i(n)$ is
the logarithm of the modulus of the eigenvalues of $\rho(n)$
corresponding to the \emph{Lyapunov splitting}
$\R^d=E_1\oplus\dots\oplus E_l$. We shall work with the natural
extension of the Lyapunov exponents to $\R^k$, that is,
$\chi:\R^k\to\R$ is a linear functional that coincides with the
Lyapunov exponent on $\Z^k$. So that for every Lyapunov exponent
$\chi=\chi_i$ we have the \emph{Lyapunov space} $E_{\chi}=E_i$ where
the eigenvalues of $\rho|E_{\chi}$ have modulus the exponential of
$\chi$. Given a Lyapunov space we denote $\chi_E$ the Lyapunov
exponent associated with $E$. Given a Lyapunov space $E_{\chi}$, we
define the \emph{complementary Lyapunov space} to be the invariant
space $\hat{E}_{\chi}$ complementary to $E_{\chi}$, i.e.
$\hat{E}_{\chi}$ is the sum of all the other Lyapunov spaces,
$E_{\chi}\oplus\hat{E}_{\chi}=\R^d$. It may be the case that two
Lyapunov exponents be positively proportional, so, given a Lyapunov
exponent $\chi$, we define the \emph{coarse Lyapunov space}
$E^{\chi}=\bigoplus_{\lambda}E_{\lambda}$ where the sum ranges over
all positive multiples $\lambda=c\chi$ of $\chi$. We define the
\emph{complementary coarse Lyapunov space} to be the invariant space
$\hat{E}^{\chi}$ complementary to $E^{\chi}$,
$E^{\chi}\oplus\hat{E}^{\chi}=\R^d$. So we shall also have the
\emph{coarse Lyapunov splitting} $\R^d=E^1\oplus\dots\oplus E^{l'}$
where each $E^i$, $1\leq i\leq l'$ is a coarse Lyapunov space. We
will call the planes $\ker\chi$ the \emph{Weyl chamber walls} and
each connected component of the complement
$\R^k\setminus\bigcup_{\chi}\ker\chi$ a \emph{Weyl chamber}. Observe
that a Weyl chamber is a cone $C\subset\R^k$ where the Lyapunov
exponents do not change sign, i.e. if $n_1, n_2\in C\cap\Z^k$ then
for every Lyapunov exponent $\chi$, $\chi(n_1)>0$ if and only if
$\chi(n_2)>0$. Given a Weyl chamber $C$, let us define the stable
space associated to any $n\in C$,
$E^s_C=\bigoplus_{\chi<0}E_{\chi}$, where the sum range over all
Lyapunov exponents that are negative on $C$. Similarly we define the
unstable space $E^u_C=\bigoplus_{\chi>0}E_{\chi}$.

See \cite{kak1} for more detailed definitions.

\end{subsubsection}
\end{subsection}

\begin{subsection}{Main theorem}
Let $\rho_*:\Z^k\to\GNZ$ be an embedding and let us denote also with
$\rho_*$ the associated standard action on $\T^N$. We shall assume
on $\rho_*$ the following properties
\begin{enumerate}[i)]
\item \label{hyp1} the coarse Lyapunov splitting coincides with the splitting of
$\R^N$ into the eigenspaces for $\rho_*$ and the eigenvalues for
$\rho_*$ are simple, in particular the coarse Lyapunov spaces
coincide with the Lyapunov spaces;
\item \label{hyp2} on each eigenspace the set of
eigenvalues form a dense subset of $\R^+$ or $\C$ depending on
wether it correspond to real or complex eigenvalues;
\item \label{hyp3} for every Weyl chamber $C$, and for every Lyapunov space
$E\subset E^s_C$ there exists an element $m\in\Z^k$ with
$\chi_E(m)<0$ and $\chi_F(m)>0$ for all other Lyapunov spaces
$F\subset E^s_C$;
\item \label{hyp4} for every Weyl chamber $C$ we want an element $m\in\Z^k$
with the following bunching property
$\chi_C^s(m)+\chi_C^{u,+}(m)-\chi_C^{u,-}(m)<0$ where $\chi_C^s(m)$
is the biggest Lyapunov exponent of $\rho_*(m)|E^s_C$ and
$\chi_C^{u,+}(m)$, $\chi_C^{u,-}(m)$ are the biggest and smallest
Lyapunov exponents of $\rho_*(m)|E^u_C$ respectively;
\end{enumerate}
\begin{theorem}\label{main}
Every linear action $\rho_*$ on the torus with the above properties
is globally rigid, that is, any Anosov action on $\T^N$ that induces
the action $\rho_*$ in homology is smoothly conjugated to it.
\end{theorem}
It is not hard to see that the actions on theorem~\ref{cent} satisfy
hypothesis \ref{hyp1})--\ref{hyp4}). Besides, if for example
$\rho_1:\Z^{k_1}\to\GZ{N_1}$ and $\rho_2:\Z^{k_2}\to\GZ{N_2}$ are
actions as in theorem~\ref{cent} then the product action
$\rho_*:\Z^{k_1+k_2}\to\GZ{N_1+N_2}$ given by
$\rho_*(n_1,n_2)=(\rho_1(n_1), \rho_2(n_2))$ also satisfies these
hypothesis and hence we can apply theorem~\ref{main} and this
product action is globally rigid.  We think that in fact if any two
actions satisfy hypothesis \ref{hyp1})--\ref{hyp4}) then their
product should also satisfy these hypothesis, the first $3$ are
easily seen, but we were not able to see how to get hypothesis
\ref{hyp4}).

Other types of actions satisfying the above hypothesis are the
following. Recall that $\mbox{Sp}(n,\Z)$ is the group of symplectic
$n\times n$ matrices with integral entries (clearly $n$ is even).
\begin{theorem}\label{centsimp}
Let $A\in\mbox{Sp}(N,\Z)$, $N\geq 4$ be a matrix with characteristic
polynomial irreducible over $\Z$, if $N=4$ assume also that $A$ has
at least one real eigenvalue. Then the standard action associated to
any finite index subgroup of $Z(A)\cap \mbox{Sp}(N,\Z)$ is globally
rigid.
\end{theorem}

Finally, the smoothness required in theorem~\ref{main} is $C^2$
although a $C^{1+\alpha}$ hypothesis would be enough, for some
$0<\alpha<1$ that a priori may depend on the action $\rho_*$.
\end{subsection}


\end{section}

\begin{section}{Smoothness of holonomies}\label{varios}
\begin{proposition}\label{conjlyap}
Let $f,g:M\to M$ be $C^{1+\mbox{\tiny{H\"older}}}$ diffeomorphisms.
Let $\mu$ be an invariant measure for $f$ and assume there is a
H\"older continuous homeomorphism $h:U\to V$ from a neighborhood $U$
of the support of $\mu$ onto $V\subset N$ such that $h\circ f=g\circ
h$. Let us call $\nu=h_*\mu$. Then, for $\mu-$a.e. $x$, the number
of negative Lyapunov exponents at $x$ are less than or equal to the
ones at $h(x)$.
\end{proposition}

In the proof we shall use the strong stable (unstable) manifold
theorem, see for instance \cite{ps}
\begin{theorem}{\bf Pesin strong stable manifold theorem.}
Let $f:M\to M$ be a $C^{1+\mbox{\tiny{H\"older}}}$ diffeomorphisms.
Let $\mu$ be an invariant measure for $f$. There is a set of full
$\mu$-measure $R_{\mu}$, the $\mu-$regular points, such that if
$x\in R_{\mu}$ and $T_xM=E_1(x)\oplus E_2(x)$ where the Lyapunov
exponents corresponding to $E_1(x)$ are negative and less than the
Lyapunov exponents corresponding to $E_2(x)$, then there is a unique
manifold $W_{E_1}(x)$ tangent to $E_1(x)$ at $x$  and characterized
as the points $y$ such that $d(y,x)\leq\epsilon(x)$ for some
$\epsilon(x)>0$ and
$$
\limsup_{n\to +\infty}\frac{1}{n}\log
d\left(f^n(y),f^n(x)\right)<\inf\{\chi_2(x),0\}
$$
where $\chi_2(x)$ is the smallest Lyapunov exponent corresponding to
$E_2(x)$.
\end{theorem}

\begin{proof}[Proof of proposition~\ref{conjlyap}]
Take a $\mu-$regular point $x$ and assume that $h(x)$ is also an
$\nu-$regular point (as $\nu=h_*\mu$, this holds for $\mu-$a.e.
point). Take the splitting $T_xM=E^s(x)\oplus E^c(x)\oplus E^u(x)$
w.r.t. negative, zero and positive Lyapunov exponents and let
$W^s(x)$ be the invariant manifold tangent to $E^s(x)$ given by the
Pesin strong stable manifold theorem. Do the corresponding
counterpart at $h(x)$. As $h$ is H\"older continuous,
$h(W^s(x))\subset W^s(h(x))$. Then, using the invariance of domain
theorem we get that $\dim(W^s(x))\leq\dim(W^s(h(x)))$ and we are
done.
\end{proof}

\begin{corollary}\label{hyphold}
In the setting of proposition~\ref{conjlyap}, if $\mu$ is a
hyperbolic measure, then $\nu$ is also a hyperbolic measure with
$\dim E^s(x)=\dim E^s(h(x))$ and $\dim E^u(x)=\dim E^u(h(x))$
$\mu-$a.e. $x$.
\end{corollary}

\begin{remark}
If $\mu$ is a hyperbolic measure, it can also be proved that if
$\chi^-_\mu\leq\log\lambda<0$ (resp. $\chi^+_\mu\geq\log\sigma>0$)
for every negative (resp. positive) Lyapunov exponent, then
$\chi^-_\nu\leq\theta\log\lambda$ (resp.
$\chi^+_\nu\geq\theta\log\sigma>0$) for every negative (resp.
positive) Lyapunov exponent, where $\theta$ is a H\"older exponent
for $h$.
\end{remark}

An improved version of the next proposition, where the existence of
the sequence $b_n$ is not needed, already appeared in \cite{sc}, we
include a proof here because it is simpler in our case.

\begin{proposition}\label{uniqueerg}
Let $f:X\to X$ be a continuous map of a compact metric space. Let
$a_n:X\to\R$, $n\geq 0$ be a sequence of continuous functions such
that $a_{n+k}(x)\leq a_n(f^k(x))+a_k(x)$ for every $x\in X$,
$n,k\geq 0$ and such that there is a sequence of continuous
functions $b_n$, $n\geq 0$ satisfying $a_n(x)\leq
a_n(f^k(x))+a_k(x)+b_k(f^n(x))$ for every $x\in X$, $n,k\geq 0$. If
$$
\inf_n \frac{1}{n}\int_X a_nd\mu<0
$$
for every ergodic $f$-invariant measure, then there is $N\geq 0$
such that $a_N(x)<0$ for every $x\in X$.
\end{proposition}

\begin{proof}
For an invariant measure $\mu$, let us call $a_n(\mu)=\int_X
a_nd\mu$. We have that $a_{n+k}(\mu)\leq a_n(\mu)+a_k(\mu)$. Now, if
$\inf_n \frac{a_n(\mu)}{n}<0$ for every ergodic $f$-invariant
measure, then the same holds for every invariant measure by the
ergodic decomposition theorem and the multiplicative ergodic
theorem. Because of the compactness of the set of invariant measures
and the properties of the $a_n(\mu)$, there is $m\geq 0$ such that
$a_m(\mu)<c<0$ for every invariant measure $\mu$. This implies that
for some $n_0>0$
$$
\sum_{j=0}^{n-1}a_m(f^j(x))<cn
$$
for every $x\in X$, $n\geq n_0$. Take $N=lm$ for some $l\geq 0$ big
enough, then
$$
\sum_{h=0}^{m-1}\sum_{i=0}^{l-1}a_m(f^{im}(f^h(x)))=\sum_{j=0}^{N-1}a_m(f^j(x))<cN
$$
Thus, we have, by the properties of the $a_n$'s, that
$$
\sum_{h=0}^{m-1}a_N(f^h(x))<cN
$$
and hence
\begin{eqnarray*}
ma_N(x)&\leq&\sum_{h=0}^{m-1}a_N(f^h(x))+a_h(x)+b_h(f^N(x))\\
&\leq&cN+m(\sup_{x\in X;\,h<m}a_h+\sup_{x\in X;\,h<m}b_h)
\end{eqnarray*}
Finally, taking $l$ big enough, as $c<0$ we get the proposition.
\end{proof}

Applying proposition~\ref{uniqueerg} to the functions $a_n(x)=\log
|D_xf^n|E|$ and $b_n(x)=\log m(D_xf^n|E)$ we get the following
immediate corollaries of the above proposition. A regular $C^1$ map
is a map whose derivative is invertible at each point.
\begin{corollary}\label{cao}
Let $f:M\to M$ be a regular $C^1$ map and $\Lambda$ a compact
invariant set. Assume $f$ leaves invariant a continuous bundle $E$
over $\Lambda$. If the Lyapunov exponents of the restriction of $Df$
to $E$ are all negative (positive) for every ergodic invariant
measure, then $Df$ contracts (expands) $E$ uniformly.
\end{corollary}
Corollary \ref{cao} already appeared in \cite{ca}. The following is
a corollary of the above and Corollary~\ref{hyphold}.
\begin{corollary}\label{invhypconj}
Let $f:M\to M$ be a diffeomorphism and $g:N\to N$ be a
$C^{1+\mbox{\tiny{H\"older}}}$ diffeomorphism. Let $\Lambda$ be a
transitive hyperbolic set for $f$ and assume there is a H\"older
continuous homeomorphism $h:U\to V$ from a neighborhood $U$ of
$\Lambda$ onto $V\subset N$ such that $h\circ f=g\circ h$. Let us
assume that $g$ leaves a continuous invariant splitting
$TM=E_1\oplus E_2$ over $h(\Lambda)=\Lambda_g$ and that it coincides
with the Lyapunov (stable$\oplus$unstable) splitting for some
(necessarily) hyperbolic $g$-invariant measure. Then $\Lambda_g$ is
a hyperbolic set for $g$.
\end{corollary}
\begin{proof}
Although in Corollary~\ref{hyphold} $f$ is assumed to be
$C^{1+\mbox{\tiny{H\"older}}}$, it is not hard to see that $\Lambda$
being a hyperbolic set, this hypothesis can be removed.
\end{proof}

Similarly we have the following corollary that states that the fact
of being an expanding map is preserved by H\"older conjugacies.
\begin{corollary}
Let $f:M\to M$ be an expanding map and $g:N\to N$ be a
$C^{1+\mbox{\tiny{H\"older}}}$ regular map. Assume there is a
H\"older continuous homeomorphism $h:M\to N$ conjugating $f$ and
$g$, i.e. $h\circ f=g\circ h$. Then $g$ is an expanding map.
\end{corollary}

The following has been recently proven by Wenxiang Sun and Zhenqi
Wang after the work of Anatole Katok \cite{k1}
\begin{theorem}\cite{sw}
Let $g:M\to M$ be a $C^{1+\alpha}$ diffeomorphism and let $\mu$ be
an ergodic hyperbolic measure. Then the Lyapunov exponents of $\mu$
can be approximated by the Lyapunov exponents of periodic orbits.
\end{theorem}

The following theorem is essentially proved in \cite{hps} and
\cite{psw}. I would like to thank Keith Burns for pointing out this
theorem to me.
\begin{theorem}
Let $f:M\to M$ be a $C^k$ diffeomorphism with an invariant splitting
$TM=E_1\oplus E_2$ satisfying
$\sup_p\frac{\|D_pf|_{E_1}\|}{m(D_pf|_{E_2})}<1$. Let us assume that
$$
\sup_x\|D_xf|_{E_1}\|\frac{\|D_xf|_{E_2}\|^r}{m(D_xf|_{E_2})}<1
$$
Then there is a $C^s$ foliation tangent to $E_1$ where $s=\min\{k-1,r\}$.
\end{theorem}


Finally we have the following corollary that get smoothness of the
strong stable foliation for an Anosov diffeomorphism form periodic
data.
\begin{corollary}\label{smoothfoliation}
Let $g$ be a $C^k$ Anosov diffeomorphism, and assume it preserves a
continuous splitting $TM=E_1\oplus E_2$ (not necessarily the
hyperbolic splitting). Given a periodic point $p$, let us call
$\chi_1^+(p)$ the biggest Lyapunov exponent of the restriction of
$Df$ to $E_1$, $\chi_2^+(p)$ the biggest Lyapunov exponent of the
restriction of $Df$ to $E_2$ and $\chi_2^-(p)$ the smallest Lyapunov
exponent of the restriction of $Df$ to $E_2$. If there is a constant
$c<0$ such that $\chi_1^+(p)-\chi_2^-(p)<c<0$ and
$\chi_1^+(p)+r\chi_2^+(p)-\chi_2^-(p)<c<0$, where $r\geq 1$, for
every periodic point $p$ then there is a $C^s$ foliation tangent to
$E_1$ where $s=\min\{k-1,r\}$.
\end{corollary}
\end{section}

\begin{section}{Smooth linearization in $\R^d$.}\label{smolin}
In this section we shall prove a result about smooth linearization
of some abelian actions in $\R^d$ that fix the origin. We shall
follow the proof of Hartman in \cite{ha} of smooth linearization of
contractions.

Take $\rho:\Z^k\to\mathit{Diff}\,^2(\R^d,0)$ an action fixing the
origin. Let us assume that there is $n_0\in\Z^k$ such that
$D_0\rho(n_0)$ is a contraction. Write $T=\rho(n_0)$ then $\R^d$
splits as a direct $D_0T$-invariant sum
$\R^d=E_{\lambda_1}\oplus\dots\oplus E_{\lambda_n}$ where
$0<\lambda_1<\dots<\lambda_n<1$ and the eigenvalues of
$D_0T|E_{\lambda_i}$ have modulus $\lambda_i$. We shall assume that
for every $i=1,\dots,n$ there is a $C^2$ manifold, $W^i$ tangent to
$E_{\lambda_i}\oplus\dots\oplus E_{\lambda_n}$ and invariant by
$\rho$ in the following sense: for every $n\in\Z^k$, there is
$\epsilon>0$ such that $\rho(n)(W^i\cap B_{\epsilon}(0))\subset
W^i$.

\begin{theorem}\label{linhart}
Let $\rho:\Z^k\to\mathit{Diff}\,^2(\R^d,0)$ be an action as above,
then there is a $C^{1+\mbox{\tiny{H\"older}}}$ diffeomorphism $h$
such that $h\circ\rho=D_0\rho\circ h$.
\end{theorem}

As in \cite {ha}, the $C^2$ condition can be relaxed to a
$C^{1+\alpha}$ hypothesis for some $0<\alpha<1$ that depends on the
eigenvalues of the action $D_0\rho$. Let us point out that the
existence of the manifolds $W^i$ is non trivial at all. The
following is an example where these manifolds are not present:
$f(x,y)=(\lambda^2x,\lambda y)$ and $g(x,y)=(x+y^2,y)$, $\lambda<1$
commute, $f$ is a linear contraction, but this $\Z^2$-action is not
linearizable, and there is no invariant manifold $W^2$ tangent to
the vertical direction.

Theorem~\ref{linhart} is of a local nature, in fact the action needs
only to be a germ of action and there will be also a smooth local
linearization. Moreover, once there is a local linearization, it can
be extended globally using the contraction $T=\rho(n_0)$ in the
obvious manner.

\begin{proof}
We shall show here how to adapt the proof in \cite{ha} to our case.
So let us describe how this proof works. Hartman's proof is
essentially by induction, he assumed the coordinates associated to
$E_{\lambda_{i+1}}\oplus\dots\oplus E_{\lambda_n}$ are already
linearized, then he found an invariant manifold tangent to
$E_{\lambda_{i+1}}\oplus\dots\oplus E_{\lambda_n}$, he make a first
conjugacy sending this invariant manifold into
$E_{\lambda_{i+1}}\oplus\dots\oplus E_{\lambda_n}$ and finally he
linearize the $E_{\lambda_i}$ coordinate, without touching the
already linearized coordinates. Finally, the first induction step is
trivial by adding some dummy coordinates.

In our case, let us make first a smooth conjugacy and assume that
the manifolds $W^i$ are already the spaces
$E_{\lambda_i}\oplus\dots\oplus E_{\lambda_n}$ and hence that the
action $\rho$ preserves this spaces.

Then, we follow the proof of Hartman. Write an $N$-vector as
$(x,y,z)$ where $x$ is an $I$-vector, $y$ a $J$-vector and $z$ a
$K$-vector and $I+J+K=d$. Let $A, B, C$ be square matrices of order
$I, J, K$ and with eigenvalues $a_1, \dots, a_I$, $b_1, \dots, b_J$,
$c_1, \dots, c_K$ respectively.

\emph{Induction hypothesis.} Assume that $T$ is written as
$$
T:x^1=Ax+X(x,y,z),\;\;\;y^1=By+Y(x,y,z),\;\;\;z^1=Cz,
$$
where the eigenvalues of $A, B, C$, satisfy
$$
0<|a_1|\leq\dots\leq |a_I|<|b_1|=\dots=|b_J|<|c_1|\leq\dots\leq
|c_K|<1;
$$
and $X, Y$ satisfy
\begin{enumerate}[(1)]
\item $X, Y$ are $C^1$ and $|(X, Y)(x,y,z)|\leq L\left(|x|+|y|+|z|\right)\left(|x|+|y|\right)$;
\item $\partial_xX, \partial_yX$ and $\partial_xY, \partial_yY$ are
uniformly Lipschitz continuous w.r.t. $(x,y,z)$;
\item $\partial_zX, \partial_zY$ are uniformly Lipschitz continuous
w.r.t. $(x,y)$;
\item $\partial_zX, \partial_zY$ are uniformly H\"older continuous
w.r.t. $z$.
\end{enumerate}

Then, as in Hartman's theorem, theorem~\ref{linhart} is proven if
the following is verified.

\emph{Induction assertion.} There exists a map $R$ of the form
$$
R:u=x,\;\;\; v=y-\varphi(x,y,z),\;\;\;w=z,
$$
where $\varphi$ satisfy
\begin{enumerate}[a)]
\item $\varphi$ is $C^1$ and $|\varphi(x,y,z)|\leq L\left(|x|+|y|+|z|\right)\left(|x|+|y|\right)$;
\item $\partial_x\varphi, \partial_y\varphi$ are uniformly Lipschitz continuous w.r.t. $(x,y,z)$;
\item $\partial_z\varphi$ is uniformly Lipschitz continuous w.r.t. $(x,y)$;
\item $\partial_z\varphi$ is uniformly H\"older continuous
w.r.t. $z$.
\end{enumerate}
$R$ is such that $F=R\circ T\circ R^{-1}$ has the form
$$
F: u^1=Au+U(u,v,w),\;\;\;v^1=Bv,\;\;\;w^1=Cw,
$$
where
\begin{enumerate}[(1)]
\item $U$ is $C^1$ and $|U(u,v,w)|\leq L(|u|+|v|+|w|)|u|$;
\item $\partial_uU$ is uniformly Lipschitz continuous w.r.t. $(u,v,w)$;
\item $\partial_vU, \partial_wU$ are uniformly Lipschitz continuous
w.r.t. $u$;
\item $\partial_vU, \partial_wU$ are uniformly H\"older continuous
w.r.t. $(v,w)$.
\end{enumerate}
We put one more assertion that is
\begin{enumerate}[(5)]
\item If $G$ has the form
$$
G: u^1=A_Gu+U_G(u,v,w),\;\;\;v^1=B_Gv+V_G(u,v,w),\;\;\;w^1=C_Gw,
$$
where $|V_G(u,v,w)|\leq L(|u|+|v|+|w|)(|u|+|v|)$ and $G$ commutes
with $F=R\circ T\circ R^{-1}$ then $V_G\equiv 0$.
\end{enumerate}

Observe that a map $\eta(a,b)$ satisfies that $|\eta(a,b)|\leq
L(|a|+|b|)|a|$ if $\eta$ is $C^1$, $\partial_a\eta$ is uniformly
Lipschitz continuous w.r.t. $(a,b)$, $\partial_a\eta(0,0)=0$ and
$\eta(0,b)=0$.

The construction of the conjugacy $R$ and the proof of assertions
(1)--(4) follows exactly the lines in \cite{ha}. Let us see the
proof of assertion (5).

We have that $F\circ G=G\circ F$ implies that $V_G\circ F^n=B^nV_G$.
Let us write the first component in $F^n$ as $(F^n)_1$, that is
$$
F^n(u,v,w)=\left((F^n)_1(u,v,w),B^nv,C^nw\right).
$$
Since $|U(u,v,w)|\leq L(|u|+|v|+|w|)|u|$ we have that for any
$\lambda>|a_K|$, $\lambda^{-n}|(F^n)_1(u,v,w)|\to 0$ if $(u,v,w)$ is
close enough to $0$. So we have that
\begin{eqnarray*}
|B^nV_G(u,v,w)|&=&|V_G(F^n(u,v,w))|\\
&\leq&L(|(F^n)_1|+|B^nv|+|C^nw|)(|(F^n)_1|+|B^nv|)
\end{eqnarray*}
where the argument of $(F^n)_1$ is $(u,v,w)$. Then, taking
$|c_K|<\mu<1$ we have that
$$
|(F^n)_1|+|B^nv|+|C^nw|\leq C\mu^n
$$
and then
$$
||b_1|^{-n}B^nV_G(u,v,w)|\leq
LC\mu^n(|b_1|^{-n}|(F^n)_1|+|b_1|^{-n}|B^nv|)
$$
Since the matrix $|b_1|^{-1}B$ has all its eigenvalues of modulus
one, we get that the norm of the matrix $|b_1|^{-n}B^n$ is bounded
between $Cn^J$ and $C^{-1}n^{-J}$ and hence
$$
C^{-1}n^{-J}|V_G(u,v,w)|\leq LC\mu^n(C+Cn^J)
$$
which gives that $V_G(u,v,w)=0$ if $(u,v,w)$ is close to $0$, then
using that $F$ is a contraction we can dispense the requirement
$(u,v,w)$ is close to $0$.

So that assertion (5) says that when we linearize $T=\rho(n_0)$ we
also linearize the whole action. Indeed, the elements of our action
satisfy the requirement on $V_G$ since they preserve the spaces
$E_{\lambda_i}\oplus\dots\oplus E_{\lambda_n}$ and also the
conjugacy $R$ preserves that spaces.

\end{proof}

To apply the theorem above we will need the following proposition

\begin{lema}\label{rigexpper}
Let $\rho_*:\Z^k\to\C\setminus\{0\}$ be a linear action induced by
complex multiplication and let $\rho:\Z^k\to\Diff{\C,0}$ be another
action. Assume that there is a homeomorphism $h:\C\to\C$ such that
$h\circ\rho=\rho_*\circ h$. Assume also that the image of $\rho_*$
is dense in $\C$. Then, $D_0\rho(n)$ has complex eigenvalues for
every $n$.
\end{lema}
\begin{proof}
If for some $n_0$, $D_0\rho(n_0)$ has its eigenvalues of the same
modulus, then applying theorem~\ref{linhart} and lemma
\ref{rigcomplex} we get the desired result. So let us assume that
$D_0\rho(n_0)$ has two eigenvalues. Then by the strong stable
manifold theorem, there is a unique $\rho(n_0)$-invariant manifold
tangent to the eigenspace of the smallest eigenvalue. Then, by
uniqueness, this manifold should be invariant by the whole action.
On the other hand, it should be dense by hypothesis, this gives a
contradiction.
\end{proof}


\end{section}

\begin{section}{Rigidity for H\"older conjugacies at periodic
orbits.}\label{righol}
In this section we shall prove that if a linear action in $\R^d$ is
H\"older conjugated to a sufficiently rich linear action then the
conjugacy should split.

\begin{definition}\label{rich}
We say that the action $\rho_*:\Z^k\to GL(d,\R)$ is \emph{rich} if
there is an element $\rho_*(n_0)$ that is a contraction and for
every coarse Lyapunov space $E_i^*$, $i=1,\dots l$ there is $n_i$
such that, $\rho_*(n_i)$ is a contraction when restricted to $E_i^*$
and an expansion on the complement.
\end{definition}
\begin{theorem}\label{lineal}
Let $\rho_*,\rho:\Z^k\to GL(d,\R)$ be linear actions on $\R^d$.
Assume $\rho_*$ is a rich action and that
$h\circ\rho(n)=\rho_*(n)\circ h$ for every $n\in\Z^k$, where
$h:\R^d\to\R^d$ is a homeomorphism H\"older continuous in a
neighborhood of the origin. Then $h$ is of the form
$$
h(x_1,\dots, x_l)=\bigl(h_1(x_1), \dots, h_l(x_l)\bigr)
$$
where $x=(x_1,\dots, x_l)$ is taken w.r.t. the coarse Lyapunov
splitting for $\rho$ and $h(x)=\bigl(h_1(x_1), \dots,
h_l(x_l)\bigr)$ is taken w.r.t. the coarse Lyapunov splitting for
$\rho_*$.
\end{theorem}

Observe that it is an implicit consequence of the theorem the fact
that $\rho$ will have a splitting $E_1\oplus\dots\oplus E_l$ that
will coincide with the coarse Lyapunov splitting.


\begin{lema}\label{ti}
Let $\rho_*:\Z^k\to GL(d,\R)$ be an action such that all its
Lyapunov exponents are positively proportional. Let $\rho:\Z^k\to
GL(d,\R)$ be another action and assume that there is an
homeomorphism $h:\R^d\to\R^d$ such that $h\circ\rho=\rho_*\circ h$.
Then, the Lyapunov exponents of $\rho$ are positively proportional
to the ones of $\rho_*$.
\end{lema}
\begin{proof}
Since $h$ conjugates $\rho$ and $\rho_*$ and both are linear, for
every  $n\in\Z^k$, $\chi_i(n)<0$ if and only if $\chi(n)<0$. Since
$\chi_i$ and $\chi$ are linear functionals the result follows.
\end{proof}

Let $\rho_*:\Z^k\to\GR{d}$ be an action, $E_1$ a coarse Lyapunov
space and $E_2$ the complementary coarse Lyapunov space,
$\R^d=E_1\oplus E_2$. Assume that there is $n_0\in\Z^k$ such that
$\rho_*(n_0)$ is a contraction and that there is $n_1\in\Z^k$ such
that $\rho_*(n_1)$ is an expansion on $E_1$ and a contraction on
$E_2$.

\begin{proposition}\label{riglinear}
Let $\rho_*:\Z^k\to GL(d,\R)$ be as above and let $\rho:\Z^k\to
GL(d,\R)$ be another action. Assume that there is $h:\R^d\to\R^d$ a
homeomorphism that is H\"older continuous in a neighborhood of the
origin and such that $h\circ\rho=\rho_*\circ h$. Then we have that
$h^{-1}(E_i)$, $i=1,2$, are complementary linear subspaces,
preserved by $\rho$ and $h_1(x_1,x_2)=h_1(x_1,0)$ where
$x=(x_1,x_2)$ are coordinates with respect to $h^{-1}(E_1)\oplus
h^{-1}(E_2)$ and $h=(h_1,h_2)$ are coordinates with respect to
$E_1\oplus E_2$
\end{proposition}

Notice that in this proposition the inverse of $h$ is not required
to be H\"older continuous at all. Moreover, it seems that the
hypothesis of being a homeomorphisms could be relaxed.

To proof the proposition we shall use the following lemma
\begin{lema}\label{laC}
Let $C, \bar C\in GL(d,\R)$ leave invariant a splitting $E_1\oplus
E_2$. Assume that $C$ and $\bar C$ are contractions and that there
is $h:\R^d\to\R^d$ a map that is H\"older continuous in a
neighborhood of the origin with H\"older exponent $\beta$ such that
$h\circ \bar C=C\circ h$. Let us call $\chi_1$ the  smallest
Lyapunov exponent of $C|_{E_1}=C_1$ and $\bar\chi_2$ the biggest
Lyapunov exponent of $\bar C|_{E_2}=\bar C_2$. If
$\beta\bar\chi_2<\chi_1$ then $h_1(x_1,x_2)=h_1(x_1,0)$ where $h_1$
is the component of $h$ in $E_1$.
\end{lema}

\begin{proof}
Denote $\bar C_1=\bar C|_{E_1}$. We may assume without loss of
generality that the neighborhood where $h$ is H\"older is the unit
ball (in fact by the conjugacy property, $h$ is H\"older in all
$\R^d$).
Take $x=(x_1,x_2)$ and take $n\geq 0$ big enough such that $|\bar
C^nx|\leq 1$, then
\begin{eqnarray*}
\bigl|h_1(x_1, x_2)-h_1(x_1, 0)\bigr|&=&\Bigl|C_1^{-n}\Bigl[h_1\bigl(\bar C_1^nx_1,\bar C_2^nx_2\bigr)-h_1\bigl(\bar C_1^nx_1,0\bigr)\Bigr]\Bigr|\\
&\leq&K\bigl\|C_1^{-n}\bigr\|\Bigl|\bar C_2^nx_2\Bigr|^{\beta}\leq
K\bigl\|C_1^{-n}\bigr\|\bigl\|\bar C_2^n\bigr\|^{\beta}|x_2|^{\beta}
\end{eqnarray*}
where $K>0$ is a generic constant. The last expression tends to $0$
as $n\to +\infty$. Indeed, take $\epsilon>0$ small such that still
$V=(1-\epsilon)\beta\bar\chi_2-(1+\epsilon)\chi_1<0$. If $n$ is big
enough, then
$$
\bigl\|C_1^{-n}\bigr\|\leq
exp\bigl(n\bigl[-\chi_1(1+\epsilon)\bigr]\bigr)
$$
and also
$$
\bigl\|\bar C_2^n\bigr\|\leq
exp\bigl(n\bigl[\bar\chi_2(1-\epsilon)\bigr]\bigr)
$$
so that
$$
\bigl\|C_1^{-n}\bigr\|\bigl\|\bar
C_2^n\bigr\|^{\beta}\leq\exp\bigl(n\bigl[\beta\bar\chi_2(1-\epsilon)-\chi_1(1+\epsilon)\bigr]\bigr)=\exp(nV)
$$
and we are done.
\end{proof}

\begin{proof}[Proof of proposition~\ref{riglinear}]
Let as assume that the H\"older exponent of $h$ is $\beta$. Take
$A=\rho_*(n_0)$ and $B=\rho_*(n_1)$ and denote $A_i=A|_{E_i}$ and
$B_i=B|_{E_i}$, $i=1,2$. As $B_1$ is an expansion and $B_2$ a
contraction, $h$ is a homeomorphism and $\rho$ is linear, by the
stable manifold theorem and the uniqueness of stable and unstable
manifolds we get that $h^{-1}(E_i)$, $i=1,2$ should be linear
$\rho$-invariant subspaces. So we may assume, without loss of
generality that $\rho$ already leaves invariant the splitting
$E_1\oplus E_2$. Let us denote also $\bar A=\rho(n_0)$, $\bar
B=\rho(n_1)$, $\bar A_i=\bar A|_{E_i}$ and $\bar B_i=\bar B|_{E_i}$,
$i=1,2$. As $h$ is a homeomorphism we have that $\bar A_2$ and $\bar
B_2$ are contractions and hence there is $a>0$ such that if
$\bar\chi_2(l,m)$ is the biggest Lyapunov exponent of $\bar
A_2^l\bar B_2^m$ then $\bar\chi_2(l,m)\leq -a(l+m)$ for $l,m\geq 0$.
Call $\chi_1$ the Lyapunov exponent for $\rho_*|E_1$ such that all
the other Lyapunov exponents when restricted to $E_1$ has rate of
proportionality less that $1$. Call $\chi_1^A=\chi_1(n_0)<0$ and
$\chi_1^B=\chi_1(n_0)>0$ then the Lyapunov exponent $\chi_1$ of
$A_1^lB_1^m$ is $\chi_1(l,m)=l\chi_1^A+m\chi_1^B$. For any $l>0$,
there is $m_l\geq 0$ such that
$$
-1\leq l\frac{\chi_1^A}{\chi_1^B}+m_l<0
$$
So $-\chi_1^B\leq l\chi_1^A+m_l\chi_1^B=\chi_1(l,m_l)<0$ and hence,
as all the other Lyapunov exponents of $\rho_*|E_1$ are positively
proportional we have that $C=A^lB^{m_l}$ is a contraction and by the
choice of $\chi_1$, that it is the smallest Lyapunov exponent of
$C|_{E_1}$. Hence we get that
$$
\beta\bar\chi_2-\chi_1=\beta\bar\chi_2(l,m_l)-\chi_1(l,m_l)\leq
-\beta a(l+m_l)+\chi_1^B
$$
Taking $l$ big enough the right hand side is negative. By lemma
\ref{laC} we get the desired property.
\end{proof}
The proof of theorem~\ref{lineal} is an immediate application of
proposition~\ref{riglinear}.

\begin{lema}
Let $\rho_*,\rho:\Z^k\to\R^+$ be linear actions on the line, and
assume that there is a continuous non constant map $h:\R^+\to\R^+$,
continuous at $0$, such that $h\circ\rho=\rho_*\circ h$. Assume also
that the image of $\rho_*$ is dense in $\R^+$. Then, there are $t>0$
and $\alpha\in\R^+$ such that $h(x)=\alpha x^t$ and $\rho_*=\rho^t$.
Moreover, if $h$ is absolutely continuous with non-zero jacobian at
$0$ then $t=1$ and hence $\rho_*=\rho$.
\end{lema}
\begin{proof}
First of all, either the image of $\rho$ is dense or discrete. If it
where discrete, then we will have a vector $v\in\Z^N$ and
$\lambda>0$ such that $\rho(n)=\lambda^{v\cdot n}$. On the other
hnd, as the image of $\rho_*$ is dense, there should be $n\in\Z^N$
such that $\rho_*(n)\neq 1$ and $v\cdot n=0$. Hence we get that for
any $x$, and such $n$, $h(x)=h(\rho(n)x)=\rho_*(n)h(x)$ which is
possible only if $h(x)=0$ and hence $h$ is trivial in which case the
proposition is trivial also. Let us assume that the image of $\rho$
is also dense. Then we have that $h(x)>0$ for every $x\neq 0$.
Hence, if $\rho_*(n)=1$ then $\rho(n)=1$, if $\rho_*(n)>1$ then
$\rho(n)>1$ and if $\rho_*(n)<1$ then $\rho(n)<1$. This is only
possible if there is $t$ positive such that $\rho_*=\rho^t$. Hence,
if we put $\alpha=h(1)$ then we get the result.
\end{proof}
An immediate corollary is the following.
\begin{corollary}\label{rigline}
Let $\rho_*,\rho:\Z^k\to\R^+$ be linear actions on the line, and
assume that there is a continuous map $h:\R\to\R$ such that
$h\circ\rho=\rho_*\circ h$. Assume also that the image of $\rho_*$
is dense in $\R^+$. Then, there are $t\geq 0$ and
$\alpha_{\pm}\in\R$ such that $h(x)=\alpha_{\pm}|x|^t$ for
$x\in\R^{\pm}$. Moreover, if $h$ is not trivial then $t\neq 0$, if
$h$ is absolutely continuous with non-zero jacobian at $0$ then
$t=1$ and hence $\rho_*=\rho$ and if the jacobian is continuous at
$0$ then $\alpha_+=-\alpha_-=\alpha$ and if $h$ preserves
orientation then $\alpha>0$.
\end{corollary}

\begin{lema}\label{rigcomplex}
Let $\rho_*:\Z^k\to\C\setminus\{0\}$ be a linear action induced by
complex multiplication and let $\rho:\Z^k\to GL(2,\R)$ be another
action. Assume that there is a homeomorphism $h:\C\to\C$ such that
$h\circ\rho=\rho_*\circ h$. Assume also that the image of $\rho_*$
is dense in $\C$. Then, after a linear conjugacy, $\rho$ is induced
by complex multiplication. Moreover, there are $t>0$,
$\alpha\in\C\setminus\{0\}$ and $a\in\R$ such that if $h$ preserves
orientation, then $h(z)=\alpha z|z|^{t-1}\exp(ia\log|z|)$ and if $h$
reverses orientation, then $h(z)=\alpha \bar
z|z|^{t-1}\exp(ia\log|z|)$. Moreover, if $h$ is absolutely
continuous at $0$ with non-zero jacobian then $t=1$ and hence
$|\rho_*|=|\rho|$ and if $h$ is differentiable at $0$, then $a=0$
and either $h$ preserves orientation and $\rho_*=\rho$ or $h$
reverses orientation and $\rho_*=\bar \rho$.
\end{lema}
\begin{proof}
Let us prove first that $\rho$ is induced by complex multiplication.
Assume by contradiction that $\rho$ leaves invariant the line
$xe^{i\theta}$, $x\in\R$. Then, using corollary~\ref{rigline} we
have that $|h(xe^{i\theta})|=\alpha x^t$ for some $\alpha>0$ and for
every $x>0$. On the other hand, as the image of $\rho_*$ is dense,
we have that the image of $xe^{i\theta}$, $x>0$ by $h$ must be dense
in $\C$. But if we take $x_n$ such that $h(x_ne^{i\theta})\to z$,
then this implies that $x_n\to\frac{|z|^{1/t}}{\alpha}$ but then it
would be impossible to approach any other point with the same
modulus of $z$. So that we may assume that $\rho$ is induced by
complex multiplication.

As $h$ is a homeomorphism we have that the image of $\rho$ is dense
in $\C$. Dividing by $h(1)$ if necessary we may assume that
$h(1)=1$. Recall that $\C$ is the universal cover of
$\C\setminus\{0\}$ with the exponential being the covering map.
Thus, we may take a lift $H$ of $h$ such that $h(e^z)=e^{H(z)}$,
$H(z+2\pi i)=H(z)+2\pi i$ and lifts $\hat{\rho}_*$ and $\hat{\rho}$
of $\rho_*$ and $\rho$ respectively in such a way that
$\hat{\rho}_*, \hat{\rho}:\Z^k\to\C$ be homomorphisms acting on $\C$
by translation, $e^{\hat{\rho}_*}=\rho_*$ and $e^{\hat{\rho}}=\rho$.
Hence $H(z+\hat{\rho}(n))=H(z)+\hat{\rho}_*(n)$. Thus our hypothesis
gives us that $\{\hat{\rho}_*(n)+l2\pi i:n\in\Z^k\,;\, l\in\Z\}$ and
$\{\hat{\rho}(n)+l2\pi i:n\in\Z^k\,;\, l\in\Z\}$ are dense. But then
$H$ should be affine, because we may see $H$ as a conjugacy between
actions by dense translations on tori. Once $H$ is affine, as
$H(0)=0$ we have that $H$ is linear and as $H(2\pi i)=2\pi i$ we get
that $H(x+yi)=tx+(ax+y)i$ for some real numbers $t\neq 0$ and $a$.
Thus we get the lemma.


\end{proof}
It seems likely that the hypothesis of $h$ being an homeomorphisms
could be relaxes.
\end{section}

\begin{section}{Putting all together.}
First let us put a corollary of sections~\ref{smolin} and
\ref{righol}. Let $\rho_*:\Z^k\to\GR{d}$ be a linear action by
semi-simple matrices. Assume that the coarse Lyapunov splitting
coincides with the splitting into eigenspaces so that each coarse
Lyapunov space has dimension one or two depending on wether it
corresponds to a real eigenvalue or a complex eigenvalue. Let us
assume also that $\rho_*$ is a rich action as in definition
\ref{rich} and that on each Lyapunov direction the set of
eigenvalues form a dense subset of $\R^+$ or $\C$ depending on
wether it correspond to real or complex eigenvaules.
\begin{theorem}\label{glorigper}
Let $\rho_*:\Z^k\to\GR{d}$ be an action as above and let
$\rho:\Z^k\to\mathit{Diff}\,^2(\R^d,0)$ be an action fixing the
origin. Assume there is a conjugacy $h:\R^d\to\R^d$,
$h\circ\rho=\rho_*\circ h$ that is H\"older continuous in a
neighborhood of $0$. Then $h$ is a diffeomorphism outside the
preimage by $h$ of the union of the $\rho_*$-complementary coarse
Lyapunov spaces. Moreover,
\begin{enumerate}[a)]
\item the corresponding Lyapunov exponents for $\rho$ and for
$\rho_*$ are proportional;
\item \label{abscont} $h$ is absolutely continuous at $0$ with
nonzero jacobian if and only if the corresponding Lyapunov exponents
for $\rho$ and for $\rho_*$ coincide and in this case $h$ is
absolutely continuous with nonzero jacobian at every point;
\item \label{diff} $h$ is differentiable at $0$ with nonzero jacobian
if and only if the corresponding eigenvalues for $\rho$ and for
$\rho_*$ coincide and in this case $h$ is a diffeomorphism.
\end{enumerate}
\end{theorem}
\begin{proof}
From the rich property for $\rho_*$, the existence of the conjugacy
$h$ and lemma~\ref{rigexpper} it follows that the action $\rho$ is
in the hypothesis of theorem~\ref{linhart}. Hence $\rho$ is smoothly
conjugated to its linear part. Then theorem~\ref{lineal} give us
that the conjugacy $h$ splits w.r.t. the coarse Lyapunov splitting.
Finally, corollary~\ref{rigline} and lemma~\ref{rigcomplex} give us
the theorem.
\end{proof}

\begin{proof}{\bf of Main Theorem~\ref{main}}
Take now a linear action $\rho_*:\Z^k\to\GNZ$ as in theorem
\ref{main} and an Anosov action $\rho$ whose action in homology is
$\rho_*$. As we said we have a H\"older continuous conjugacy $h$
such that $h\circ\rho=\rho_*\circ h$ and we want to prove that it is
smooth. Let $n_0$ be such that $f=\rho(n_0)$ is an Anosov
diffeomorphism, let $A=\rho_*(n_0)$ and $C$ be the $\rho_*$-Weyl
chamber such that $n_0\in C$. Let $p$ be a periodic point for $f$
and $\Gamma_p$ the stabilizer of $p$, that is, the finite index
subgroup of $\Z^k$ that leave $p$ fixed. Let $W^s(p)$ be the stable
manifold of $p$ for $f$, we can identify it with $\R^d$ and we have
that the restriction of $\rho$ to $\Gamma_p$ leave $W^s(p)$
invariant. We have that $h(p)$ is a periodic point for $\rho_*$ and
that $\Gamma_p=\Gamma_{h(p)}$ so we can work also with the
restriction of $\rho_*$ to $\Gamma_p$. The stable space of $h(p)$
for $A$, $h(p)+E^s_C$ can also be identified with $\R^d$. Hence we
can work with the actions $\rho$ and $\rho_*$ induced on this $\R^d$
and fixing the origin. They are H\"older conjugated by the
restriction of $h$ to these stable manifolds. We can see that
hypothesis~\ref{hyp1}),~\ref{hyp2}) and~\ref{hyp3}) in theorem
\ref{main} guaranty that we are in the hypothesis of theorem
\ref{glorigper}. So we have that $h$ restricted to the stable
manifold of $p$ is a diffeomorphism outside the complementary
Lyapunov spaces. Let us see that it is in fact absolutely continuous
and hence that the Lyapunov exponent at $p$ of the restriction of
$\rho$ to $W^s(p)$ coincide with the Lyapunov exponents of the
restriction of $\rho_*$ to $E^s_C$.

Let us work in the universal covering $\R^N$. We shall use the same
notation for the objects in the torus or in $\R^N$ whenever this
leads not to confusion. It is not hard to see that there is
$n\in\Z^N$ such that $E^u_C\cap (L+n)=\emptyset$ for every
complementary Lyapunov space $L\subset E^s_C$. Let us define the
holonomy map $\pi^*_n:h(p)+E^s_C\to h(p)+n+E^s_C$ sliding along
$E^u_C$, that is
$\pi^*_n(y)=\left(h(p)+n+E^s_C\right)\cap\left(y+E^u_C\right)$,
$\pi^*_n$ is an affine map and hence it is smooth. The choice of $n$
implies that $\pi^*_n(h(p))\notin L$ for every complementary
Lyapunov space $L\subset E^s_C$. We have also the unstable holonomy
for $f$, $\pi_n:W^s(p)\to W^s(p+n)=W^s(p)+n$ defined by
$\pi_n(x)=\left(W^s(p)+n\right)\cap W^u(x)$. Since $f$ is
$C^{1+\mbox{\tiny{H\"older}}}$ it follows that $\pi_n$ is absolutely
continuous with nonzero jacobian. On the other hand, as the
conjugacy send the stable and unstable foliations into the
corresponding ones and it is homotopic to the identity we have that
it conjugates the unstable holonomies, that is
$h\circ\pi_n=\pi^*_n\circ h$ and hence $\pi_n(p)$ is not in the
preimage of the complementary Lyapunov spaces. So we have that
$h=(\pi^*_n)^{-1}\circ h\circ\pi_n$ and hence $h$ is written in a
neighborhood of $p$ as the composition of $\pi_n$ that is absolutely
continuous, $h$ restricted to a neighborhood of $\pi_n(p)$ that is a
diffeomorphism by the choice of $n$ and the inverse of $\pi^*_n$
that is also a diffeomorphism so we have that in a neighborhood of
$p$, $h$ is the composition of $\pi_n$ with a diffeomorphism and
hence it is absolutely continuous in a neighborhood of $p$ and by
conclusion~\ref{abscont}) of theorem~\ref{glorigper} we have that
$h$ is absolutely continuous when restricted to $W^s(p)$ and the
Lyapunov exponents for $\rho$ restricted to $W^s(p)$ at $p$ coincide
with the ones of $\rho^*$ restricted to $E^s_C$. Working with
$f^{-1}$ we get the same for the unstable manifold and as $p$ was
arbitrary we have that the Lyapunov exponents at any point for
$\rho$ coincide with the ones of $\rho_*$.

On the other hand, using corollary~\ref{invhypconj} we get that
$\rho(n)$ is Anosov for any element $n\in C$, the $\rho^*$-Weyl
chamber containing $n_0$. Take an element $n_1\in C$ satisfying the
hypothesis~\ref{hyp4}) of theorem~\ref{main}, by corollary
\ref{smoothfoliation} we have that the stable foliation is smooth,
similarly working with $f^{-1}$ and $-n_0\in -C$ we get that the
unstable foliation is smooth. Hence the holonomy $\pi_n$ is smooth
and hence by the same argument we used to prove that $h$ was
absolutely continuous but now using conclusion~\ref{diff}) of
theorem~\ref{glorigper} we get that $h$ restricted to $W^s(p)$ is
smooth at $p$ and hence $h|W^s(p)$ is a diffeomorphism. Similarly
$h|W^u(p)$ is a diffeomorphism. Finally as the stable and unstable
foliations are smooth we get that $h$ is a diffeomorphism and we are
done.
\end{proof}

We want to mention also that from a careful reading of the proof, it
can be seen that the $C^1$ distance of the conjugacy to the identity
depends only on the $C^2$ distance of the Anosov element of the
action to the linear one and some properties on the linear action.
In fact, this estimate comes from theorem~\ref{linhart} and from the
bounds on the regularity of the invariant foliations that are
controlled since the eigenvalues of the Anosov element coincide with
the linear ones.

\end{section}

\bibliographystyle{alpha}

\begin{thebibliography}{99}

\bibitem{br}
M.I. Brin, \emph{Nonwandering points of Anosov diffeomorphisms.}
Ast\'erisque {\bf 49} (1977), 11--18.


\bibitem[Ca]{ca}
Y. Cao, \emph{Non-zero Lyapunov exponents and uniform
hyperbolicity.} Nonlinearity 16 (2003) 1473-–1479.


\bibitem{fr}
J. Franks, \emph{Anosov diffeomorphisms.} 1970 Global Analysis Proc.
Sympos. Pure Math, Vol XIV, Berkeley, Calif. (1968), 61--93.

\bibitem{fu}
H. Furstenberg, \emph{Disjointness in ergodic theory, minimal sets,
and a problem in diophantine analysis.} Math. Syst. Theory, {\bf 1}
(1967), 1--49.

\bibitem{ha}
P. Hartman, \emph{On local homeomorphisms of Euclidean subspaces.}
Bol. Soc. Mat. Mexicana {\bf 5} (1960), 220--241.




\bibitem{hps}
M. Hirsch, C. Pugh and M. Shub, \emph{Invariant manifolds.} Lecture
Notes in Mathematics, Vol. 583. Springer-Verlag, Berlin-New York,
1977.


\bibitem{jo}
A. Johnson, \emph{Measures on the circle invariant under
multiplication by a nonlacunary subsemi- group of the integers.}
Israel J. of Math. {\bf 77} (1992), 211--240.

\bibitem{kak1}
B. Kalinin and A. Katok, \emph{Invariant measures for
actions of higher rank abelian groups.} Proc. Symp. Pure Math, {\bf
69}, (2001), 593--637.

\bibitem{kak}
B. Kalinin and A. Katok, \emph{Measure rigidity beyond uniform
hyperbolicity: Invariant Measures for Cartan actions on Tori.}
Preprint (2006).

\bibitem{kakrh}
B. Kalinin, A. Katok and F. Rodriguez Hertz, \emph{Nonuniform
measure rigidity.} in preparation.

\bibitem{kas}
B. Kalinin and R. Spatzier, \emph{On the classification of Cartan
actions.} to appear in GAFA, Geometric And Functional Analysis.

\bibitem{k1}
A. Katok, \emph{Lyapunov exponents, entropy and periodic orbits for
diffeomorphisms.} Inst. Hautes Études Sci. Publ. Math. No. {\bf 51}
(1980), 137--173.


\bibitem{kl}
A. B. Katok and J. W. Lewis, \emph{Local rigidity for certain groups
of toral automorphism.} Isr. J. of Math. {\bf 75} (1991), 203--241.

\bibitem{klz}
A. B. Katok, J. W. Lewis and R. J. Zimmer, \emph{Cocycle
superrigidity and rigidity for lattice actions on tori.} Topology
{\bf 35} (1996), 27--38.

\bibitem{krh}
A. Katok and F. Rodriguez Hertz, \emph{Uniquness of large invariant
measures for $\Z^k$ actions with Cartan homotopy data.} Preprint
(2006).

\bibitem{ks}
A. Katok and R. Spatzier, \emph{Differential rigidity of
Anosov actions of higher rank abelian groups and algebraic lattice
actions.} Tr. Mat. Inst. Steklova {\bf 216} (1997), Din. Sist. i
Smezhnye Vopr., 292–319; translation in Proc. Steklov Inst. Math.
1997, no. 1 {\bf 216}, 287--314

\bibitem{ks2}
A. Katok and R. J. Spatzier, \emph{Invariant Measures for Higher
Rank Hyperbolic Abelian Actions.} Ergod. Th. \& Dynam. Syst. {\bf
16} (1996), 751–-778.

\bibitem{ma}
A. Manning, \emph{There are no new Anosov diffeomorphisms on tori.}
Amer. J. Math. {\bf 96} (1974), 422--429.

\bibitem{ne}
S. Newhouse, \emph{On codimension one Anosov diffeomorphisms.} Amer.
J. Math. {\bf92} (1970), 761--770.

\bibitem{ps}
C. Pugh and M. Shub, \emph{Ergodic attractors.} Trans. Amer. Math.
Soc. {\bf 312} (1989), no. 1, 1--54.

\bibitem{psw}
C. Pugh, M. Shub and A. Wilkinson, \emph{H\"older Foliations.} Duke
Math. J. {\bf 86} (1997), no. 3, 517--546.

\bibitem{rh}
F. Rodriguez Hertz, \emph{Global rigidity of $\Z^2$ Cartan Actions
on $\T^3$.} Preprint (2001).

\bibitem{ru}
D. Rudolph, \emph{$\times 2$ and $\times 3$ invariant measures and
entropy.} Ergod. Th. \& Dynam. Syst. {\bf 10} (1990), 395--406.

\bibitem[Sc]{sc}
S.J. Schreiber, \emph{On growth rates of sub-additive functions for
semi-flows.} J. Differential Equations 148 (1998) 334–-350.

\bibitem{sw}
W. Sun and Z. Wang, \emph{Lyapunov exponents of hyperbolic measures
and hyperbolic periodic orbits.} Preprint (2005)

\end{thebibliography}

\end{document}